\documentclass{ws-ijm}

\usepackage [latin1]{inputenc}
\usepackage[all]{xy}
\usepackage [english]{babel}

\usepackage{amsmath}
\usepackage{amsfonts}
\usepackage{amssymb}
\usepackage{graphics}
\usepackage{graphicx}

\newtheorem{teorema}{Theorem}[section]
\newtheorem{propos}[teorema]{Proposition}
\newtheorem{corol}[teorema]{Corollary}
\newtheorem{ex}{Example}[section]
\newtheorem{rem}{Remark}[section]
\newtheorem{defin}[teorema]{Definition}

\def\bt{\begin{teorema}}
\def\et{\end{teorema}}
\def\bp{\begin{propos}}
\def\ep{\end{propos}}
\def\bl{\begin{lemma}}
\def\el{\end{lemma}}
\def\bc{\begin{corol}}
\def\ec{\end{corol}}
\def\br{\begin{rem}\rm}
\def\er{\end{rem}}
\def\bex{\begin{ex}\rm}
\def\eex{\end{ex}}
\def\bd{\begin{defin}}
\def\ed{\end{defin}}

\def\R{{\mathbb R}}   \def\a {\alpha} 
\def\N{{\mathbb N}}      \def\e{\varepsilon}
\def\C{{\mathbb C}}

\def\P{{\mathbb P}}
\def\der{\partial/\partial} \def\oli{\overline}

\def\O{\Omega}
\def\S{\Sigma}

\begin {document}
\markboth{Giuseppe Della Sala, Alberto Saracco}
{Non compact boundaries of complex analytic varieties}

\catchline{}{}{}{}{}

\title{NON COMPACT BOUNDARIES OF COMPLEX ANALYTIC VARIETIES}
\author{GIUSEPPE DELLA SALA, ALBERTO SARACCO}
               \address{Scuola Normale Superiore, Piazza dei Cavalieri, 7 - I-56126 Pisa, Italy, \\ e-mail g.dellasala@sns.it, a.saracco@sns.it}
               \maketitle
               
               \begin{abstract}We treat the boundary problem for complex varieties with isolated singularities, of dimension greater than one, which are
               contained in a certain class of strongly pseudoconvex, not necessarily bounded open subsets of $\C^n$. We deal with the problem by cutting with a
family of complex hyperplanes and applying the classical Harvey-Lawson's theorem for the bounded case \cite{HL}.
               \end{abstract}
               

\section {Introduction}
Let $M$ be a smooth and oriented $(2m+1)$-real submanifold of some
$n$-complex manifold $X$. A natural question
arises, whether $M$ is the boundary of an ($m+1$)-complex analytic
subvariety of $X$. This problem, the so-called \textit{boundary
problem}, has been widely treated over the past fifty years when
$M$ is compact and $X$ is $\C^n$ or $\C\P^n$.

The case when $M$ is a compact, connected curve in $X=\C^n$
($m=0$), has been first solved by Wermer~\cite{W} in 1958. Later
on, in 1975, Harvey and Lawson in~\cite{HL} and~\cite{HL3} solved
the boundary problem in $\C^n$ and then in
$\C\P^n\setminus\C\P^r$, in terms of holomorphic chains, for any $m$. The
boundary problem in $\C\P^n$ was studied by Dolbeault and Henkin,
in~\cite{DH} for $m=0$ and in~\cite{DH2} for any $m$. Moreover, in
these two papers the boundary problem is dealt with also for
closed submanifolds (with negligible singularities) contained in
$q$-concave (i.e.\ union of $\C\P^q$'s) open subsets of $\C\P^n$.
This allows $M$ to be non compact. The results in~\cite{DH}
and~\cite{DH2} were extended by Dinh in~\cite{D99}.\vspace{0,3cm}

The main theorem proved by Harvey and Lawson in~\cite{HL} is that
if $M\subset\C^n$ is compact and maximally complex then $M$ is the
boundary of a unique holomorphic chain of finite
mass~\cite[Theorem 8.1]{HL}. Moreover, if $M$ is contained in the
boundary $b\Omega$ of a strictly pseudoconvex domain $\Omega$ then
$M$ is the boundary of a complex analytic subvariety of $\Omega$,
with isolated singularities~\cite{HL2} (see also~\cite{G}). The
aim of this paper is to generalize this last result to a non
compact, connected, closed and maximally complex submanifold $M$ of
the connected boundary $b\Omega$ of an unbounded weakly
pseudoconvex domain $\Omega\subset \C^n$. The pseudoconvexity of $\Omega$ is needed both for the local result and to prove that the singularities are isolated.

Maximal complexity of $M$ and extension theorem for $CR$ functions (see \cite{Le}) allow us to prove the following semi global
result (see Corollary~\ref{L4}). Assume that $n\geq 3$, $m\geq 1$
and the Levi form $\mathcal L(b\O)$ of $b\O$ has at least $n-m$
positive eigenvalues at every point $p\in M$. Then
\begin{enumerate}\item[] there exist
a tubular open neighborhood $I$ of $b\Omega$ and a
complex submanifold $W_0$ of $\oli \O \cap I$ with boundary, such
that $bW_0\cap b\Omega=M$, i.e. a complex manifold $W_0 \subset
I\cap \O$ such that the closure $\oli W_0$ of $W_0$ in $I$ is a
smooth submanifold with boundary $M$.\end{enumerate}

A very simple example (see Example~\ref{E1}) shows that in general the semi global result
fails to be true for $m=0$.

In order to prove that $W_0$ extends to a complex analytic
subvariety $W$ of $\Omega$ with boundary $M$ we first treat the
case when $\Omega$ is convex and does not contain straight lines.
This is the crucial step. For technical reasons we divide the
proof in two cases: $m\geq2$ and $m=1$. We cut $\overline\Omega$
by a family of real hyperplanes $H_\lambda$ which intersect $M$
along smooth compact submanifolds. Then the natural foliation on
each $H_\lambda$ by complex hyperplanes induces on $M\cap
H_\lambda$ a foliation by compact maximally complex ($2m-1$)-real
manifolds $M'$. Thus a natural way to proceed is to apply
Harvey-Lawson's theorem to each $M'$ and to show that the family
$\{W'\}$ of the corresponding Harvey-Lawson solutions actually
organizes in a complex analytic subvariety $W$, giving the desired
extension (see Theorem~\ref{MT}). This is done by following an idea of Zaitsev (see Lemma \ref{Zaitsev}).

The same method of proof is used in the last section in order to
treat the problem when $\Omega$ is pseudoconvex. In this case,
$M$ is requested to fulfill an additional condition.
Precisely,
\begin{enumerate}\item[($\star $)]  if $\oli
M^\infty$ denotes the closure of $M \subset \C^n\subset \C\P^n$ in
$\C\P^n$, then there exists an algebraic hypersurface $V$ such
that $V\cap\oli M^\infty=\emptyset$. \end{enumerate} Equivalently
\begin{enumerate}\item[($\star'$)] there exists a polynomial $P\in \C[z_1,\ldots,z_n]$ such that $$M \subset\left\{z\in \C^n:\vert P(z)\vert^2>(1+\vert z\vert^2)^{{\rm deg} P}\right\}.$$\end{enumerate}
A similar condition was first pointed out by Lupacciolu~\cite{L2}
in studying the extension problem for $CR$ functions in unbounded
domains. It allows us to build a nice family of hypersurfaces,
which play the role of the hyperplanes in the convex case, and so
to prove the main theorem of the paper: \bt\label{GT} Let $\O$ be
an unbounded domain in $\C^n$ $(n\geq 3)$ with smooth boundary
$b\O$ and $M$ be a maximally complex closed $(2m+1)$-real
submanifold $(m \geq 1)$ of $b\O$. Assume that
\begin{enumerate}
\item [\emph{(i)}] $b\O$ is weakly pseudoconvex and the Levi form $\mathcal
L(b\O)$ has at least $n-m$ positive eigenvalues at every point of
$M$;
\item[\emph{(ii)}] $M$ satisfies condition $(\star)$.
\end{enumerate} Then there exists a unique
$(m+1)$-complex analytic subvariety $W$ of $\O$, such that $bW =
M$. Moreover the singular
locus of $W$ is discrete and the closure of $W$ in $\oli \O
\setminus {\sf Sing} \, (W)$ is a smooth submanifold with boundary
$M$. \et

We do not deal with the 1-dimensional case. There are two
different kinds of difficulties. First of all, a semi global strip as in
Corollary~\ref{L4} may not exist (see Example~\ref{E1}).
Secondly, even though it does exist, it could be non extendable to
the whole $\Omega$ (see Example~\ref{E2}) and it is not clear at all
how it is possible to generalize the \textit{moments condition}
(see~\cite{HL}).

\vspace{0.3 cm}

Another similar approach can be followed to treat the
\textit{semi-local} boundary problem, i.e. given an open subset $U$
of the boundary of $\O$, find an open subset $\O'\subset \O$ such
that, for any maximally complex submanifold $M\subset U$, there
exists a complex subvariety $W$ of $\O'$ whose boundary is $M$. We
deal with this problem in a work in preparation (see \cite{DS}).

\section{Definitions and Notations}

We briefly recall some well known notions of $CR$ geometry that
will be used in the paper.

Let $N\subset\C^n$ be a smooth connected real submanifold, and let
$p\in N$. We denote by $T_p(N)$ the tangent space of $N$ at the
point $p$, and by $H_p(N)$ the holomorphic tangent space of $N$ at
the point $p$.

A ($2k+1$)-real submanifold $N\subset\C^n$, $k\geq1$, is said to
be a \textit{$CR$ submanifold} if ${\sf dim}_\C H_p(N)$ is constant
along $N$. When this is the case, $H(N)=\cup_p H_p(N)$ is a subbundle of the tangent bundle $T(N)$. If ${\sf dim}_\C H_p(N)$ is the greatest possible, i.e.\
${\sf dim}_\C H_p(N) = k$ for every $p$, $N$ is said to be
\textit{maximally complex}.

A $C^\infty$ function $f:N\to\C$ is said to be a
\textit{$CR$ function} if for a $C^\infty$ extension (and
hence for any) $\widetilde f: U\to \C$ ($U$ being a neighborhood
of $N$) we have
\begin{equation}\label{1}\left(\oli\partial\widetilde f\right)|_{H(N)}\ =\ 0.\end{equation}
In particular the restriction of a holomorphic function to a $CR$
submanifold is a $CR$ function. It is immediately seen that $f$ is
$CR$ if and only if
\begin{equation}
df\wedge(dz_1\wedge \ldots \wedge dz_n)|_N = 0.
\end{equation}

Similarly $N$ is maximally complex if and only if 
$$
(dz_{j_1}\wedge \ldots \wedge dz_{j_{k+1}})|_N = 0,$$
for any $(j_1,\ldots,
j_{k+1})\in\left\{1,\ldots,n\right\}^{k+1}$.

Finally we observe that the boundary $M$ of a complex submanifold
$W$ with ${\sf dim}_\C W > 1$ is maximally complex. Indeed, for
any $p\in bW=M$, $T_p(bW)$ is a real hyperplane of $T_p(W)=H_p(W)$
and so is $J(T_p(bW))$. Hence $H_p(bW)=T_p(bW)\cap J(T_p(bW))$ is
of real codimension $2$ in $H_p(W)$. 

If ${\sf dim}_\C W=1$ and $bW$ is compact then for any holomorphic
$(1,0)$-form $\omega$ we have
$$\int_{M}\omega\ =\ \iint_W
d\omega\ =\ \iint_W
\partial\omega\ =\ 0,$$
since $\partial\omega|_W \equiv 0$. This
condition for $M$ is called \textit{moments condition} (see
\cite{HL}).

By the same arguments, a ($2n-1$)-real submanifold of $\C^n$ is maximally complex.

\section{The Local and Semi Global Results}\label{local}

The aim of this section is to prove the local result. Given a
smooth real hypersurface $S$ in $\C^n$, we denote by $\mathcal L_p(S)$
the Levi form of $S$ at the point $p$.  Let $0$ be a point of $M$. We
have the following inclusions of tangent spaces: $$ \C^n\
\supset\ T_0(S)\ \supset\ H_0(S)\ \supset\ H_0(M);$$ $$
\phantom{\C^n\ \supset}\ T_0(S)\ \supset\ T_0(M)\ \supset\
H_0(M).$$

\bl\label{wk} Let $M$ be a maximally complex submanifold of a
smooth real hypersurface $S$, ${\sf dim}_\R M=2m+1$, $m\geq1$, $0\in M$. Suppose
that $\mathcal L_0(S)$ has at least $n-m$ eigenvalues of the same
sign. Then $$H_0(S)\not\supset T_0(M).$$ \el
\begin{proof} Should the
thesis fail we would have the following chain of inclusions $$ \C^n
\supset T_0(S) \supset H_0(S) \supset T \supset T_0(M) \supset
H_0(M),$$ where $T$ is the smallest complex space containing
$T_0(M)$ (since $M$ is maximally complex, ${\sf dim}_\C T=m+1$). Hence, we may choose in a neighborhood of $0$ local complex coordinates $z_k=x_k
+ i y_k$, $k=1,\ldots,m+1$, $w_l=u_l + i v_l$, $l=m+2,\ldots,n$,
in such a way that:
\begin{itemize}
\item $H_0(M) = {\sf span} \ (\der x_k, \der y_k)$, $k=1,\ldots,m$
\item $T_0(M) = {\sf span} \ (\der x_k, \der y_k, \der x_{m+1})$, $k=1,\ldots,m$
\item $T = {\sf span} \ (\der x_k, \der y_k)$, $k=1,\ldots,m+1$
\item $H_0(S) = {\sf span} \ (\der x_k, \der y_k, \der u_l, \der v_l)$, $k=1,\ldots,m+1$,
$l=m+2,\ldots,n-1$, if $m+2\leq n-1$ \\ or
\item $H_0(S) = T$, if $m=n-2$;
\item $T_0(S) = {\sf span} \ (\der x_k, \der y_k, \der u_l, \der v_l, \der u_n)$, $k=1,\ldots,m+1$,
$l=m+2,\ldots,n-1$, if $m+2\leq n-1$ \\ or
\item  $T_0(S) = {\sf span} \ (\der x_k, \der y_k, \der u_n)$ $k=1,\ldots,m+1$, if
$m=n-2$.
\end{itemize}

We denote by $z$ the first $m+1$ coordinates, by $\hat z$ the
first $m$, and by $\pi$ the projection on $T$; $\pi$ is obviously
a local embedding of $M$ near $0$, and we set $M_0 = \pi(M)$.\\
Locally  at $0$, $S$ is a graph over its tangent space: $$S=\{v_n
= h(u_n,u_j,v_j,x_i,y_i)\}.$$ Observe that the Levi form of $h$
has $n-m$ eigenvalues of the same sign. In order to obtain a similar
description of $M$, we proceed as follows. First, we have $$M_0 =
\{ (\hat z, z_{m+1}): y_{m+1}=\varphi(\hat z,x_{m+1})\}.$$ Then,
we choose $f_j(\hat z,x_{m+1}) = f_j^1(\hat z,x_{m+1}) +
if_j^2(\hat z,x_{m+1})$ (where $f^1_j$ and $f^2_j$ are
real-valued) defined in a neighborhood of $M_0$ in $T$ in such a
way that $$M= \{w_{m+2} = f_{m+2}(\hat z,x_{m+1}), \ldots, w_{n} =
f_{n}(\hat z,x_{m+1}) \}.$$ Observe that the function
$(f_{m+2}(\hat z,x_{m+1}), \ldots, f_{n}(\hat z,x_{m+1}))$ is just
$\pi^{-1}|_{M_0}$, and since $M$ is maximally complex it has to be
a $CR$ map.

By hypothesis, the following equation holds in a neighborhood of
$0$: $$f_n^2(\hat z,x_{m+1}) = h\left(f_n^1(\hat z,x_{m+1}),
f_j^k(\hat z,x_{m+1}),\hat z,x_{m+1}\right).$$ After a computation
of the second derivatives, taking into account that all first
derivatives of $f_j^k$, of $h$ and of $\varphi$ vanish in the
origin, we obtain 
$$ \frac{\partial^2 f_n^2}{\partial z_j\partial
\overline z_k}(0)\ =\ \frac{\partial^2 h}{\partial z_j\partial
\overline z_k}(0),$$
i.e.\ the Levi form of $h$ and $f_n^2$ coincide
in $H_0(M)$. By hypothesis $\mathcal L_0(h)$ is strictly positive
definite on a non-zero subspace of $H_0(M)$. We shall obtain a
contradiction by showing that $\mathcal L_0(f_n)$ (and hence $\mathcal L_0(f_n^2)$) vanishes on
$H_0(M)$. Let $\xi\in H_0(M)$. We may assume (up to unitary linear
transformation of coordinates of $H_0(M)$) that $\xi =\partial /
\partial z_1$.

Set $f\doteqdot f_n$. Then, since $f$ is a $CR$ function on $M_0$, we have:
$$\frac{\partial}{\partial \oli z_k}f(\hat z, x_{m+1}) = -\a(\hat
z, x_{m+1}) \frac{\partial}{\partial \oli z_k}\varphi(\hat z,
x_{m+1}),\ \ k=1,\ldots, m $$ and $$ \frac{\partial}{\partial \oli
z_{m+1}}f(\hat z, x_{m+1}) = -i\a(\hat z, x_{m+1}) + \a(\hat z,
x_{m+1}) \frac{\partial}{\partial x_{m+1}}\varphi(\hat z,
x_{m+1}),$$ where $\a(\hat z, x_{m+1})$ is a complex valued function. Differentiating and calculating in $0$ we obtain
\begin{equation} \label{prima}
\frac{\partial^2 f}{\partial z_1 \partial \oli{z_1}}(0) = \a(0) \frac{\partial^2 \varphi}{\partial z_1 \partial \oli{z}_1}(0),
\end{equation}

\begin{equation} \label{seconda}
0 = \frac{\partial f}{\partial x_{m+1}}(0) = i\a(0),
\end{equation}
i.e. $\a(0) = 0$. From (\ref{prima}) we deduce that $\partial^2 f / \partial z_1 \partial \oli z_1 (0) = 0$. Contradiction.
\end{proof}

\bl\label{L3} Under the hypothesis of Lemma \ref{wk}, assume that $S$
is the boundary of an unbounded domain $\O\subset \C^n$, $0\in M$
and that the Levi form of $S$ has at least $n-m$ positive
eigenvalues. Then
\begin{enumerate}
\item[\emph{(i)}] there exists an open neighborhood $U$ of $0$ and an $(m+1)$-complex submanifold $W_0\subset U$
 with boundary, such that $bW_0=M\cap U$;
\item[\emph{(ii)}] $W_0\subset\O\cap U$.
\end{enumerate}
\el \begin{proof} To prove the first assertion, observe that to obtain $\mathcal L^M_0(\zeta_0,\oli \zeta_0)$
it suffices to choose a smooth local section $\zeta$ of $H_0(M)$
such that $\zeta(0) = \zeta_0$ and compute the
projection of the bracket $[\zeta,\oli\zeta](0)$ on the real part
of $T_0(M)$. By hypothesis, the intersection of the space
where $\mathcal L_0(S)$ is positive with $H_0(M)$ is non empty; take $\eta_0$ in this intersection. Then $\mathcal L_0^M(\eta_0, \oli\eta_0)\neq 0$. Suppose, by contradiction, that the bracket $[\eta,\oli\eta](0)$ lies in $H_0(M)$, i.e.\ its projection on the real part of the tangent of $M$ is zero. Then, if $\widetilde{\eta}$ is a local smooth extension of the field $\eta$ to $S$, we have $[\widetilde{\eta},\oli{\widetilde{\eta}}](0)= [\eta,\oli\eta](0)\in H_0(M)$. Since $H_0(M)\subset H_0(S)$, this would mean that the Levi form of $S$ in $0$ is zero in $\eta_0$. 
Now, we project (generically) $M$ over a $\C^{m+1}$ in such a way that the
projection $\pi$ is a local embedding near $0$: since the
restriction of $\pi$ to $M$ is a $CR$ function, and since the Levi
form of $M$ has - by the arguments stated above - at least one
positive eigenvalue, it follows that the Levi form of $\pi(M)$ has
at least one positive eigenvalue. Thus, in order to obtain $W_0$,
it is sufficient to apply the Lewy extension theorem \cite{Le} to
the $CR$ function $\pi^{-1}|_M$.

As for the second statement, we observe that the projection by
$\pi$ of the normal vector of $S$ pointing towards $\O$ lies into
the domain of $\C^{m+1}$ where the above extension $W_0$ is
defined. Indeed, the extension result in \cite{Le} gives a holomorphic function in the connected component of (a neighborhood
of $0$ in) $\C^n \setminus \pi(M)$ for which $\mathcal
L_0(\pi(M))$ has a positive eigenvalue when $\pi(M)$ is oriented
as the boundary of this component. This is precisely the component
towards which the projection of the normal vector of $S$ points
when the orientations of $S$ and $M$ are chosen accordingly. This
fact, combined with Lemma \ref{wk} (which states that any
extension of $M$ must be transverse to $S$) implies that locally
$W_0\subset\O\cap U$.
\end{proof}

\bc[Semi global existence of $W$]\label{L4} Under the same hypothesis of Lemma~\ref{L3}, there exist an open tubular neighborhood $I$ of
$S$ in $\oli \O$ and an $(m+1)$-complex submanifold $W_0$ of
\ $\oli\O \cap I$, with boundary, such that $S\cap bW_0=M$. \ec
\begin{proof} By Lemma~\ref{L3}, for each point $p\in M$, there exist a
neighborhood $U_p$ of $p$ and a complex manifold
$W_p\subset\oli\O\cap U_p$ bounded by $M$. We cover $M$ with
countable many such open sets $U_i$, and consider the union
$W_0=\cup_i W_i$. $W_0$ is contained in the union of the $U_i$'s,
hence we may restrict it to a tubular neighborhood $I_M$ of $M$.
It is easy to extend $I_M$ to a tubular neighborhood $I$ of $S$.

The fact that $W_{i}|_{U_{ij}}=W_{j}|_{U_{ij}}$ if $U_i\cap
U_j=U_{ij}\neq\emptyset$ immediately follows from the construction made in
Lemma~\ref{L3}, in view of the uniqueness of the holomorphic extension of
$CR$ functions. \end{proof}

\bex\label{E1} Corollary~\ref{L4} could be restated by saying that
if a submanifold $M\subset S$ (${\sf dim}_\R M \geq 3$) is locally
extendable at each point as a complex manifold, then (one side of) the extension
lies in $\O$. This is no longer true, in general, for curves, as
shown in $\C^n_{(z_1,\ldots,z_{n-1},w)}$, $z_k=x_k+iy_k$, $w=u+iv$, by the following case:
$$S\
=\ \left\{v=u^2+\sum_k\left|z_k\right|^2\right\}, \ \O\ =\
\left\{v > u^2+\sum_k\left|z_k\right|^2\right\},$$ $$M\ =\
\left\{y_1=0,\ v=x_1^2,\ u=0,\ z_2=\cdots=z_{n-1}=0\right\}$$ and
$$W\ =\ \left\{w=iz_1^2,\ z_2=\cdots=z_{n-1}=0\right\};$$ we have
that $S\cap W = M$ and $W \subset \C^n \setminus \O$. \eex

\br Suppose that $S$ is strongly pseudoconvex and choose, in
$\C^n_{(z_1,\ldots,z_n)}$, a local strogly plurisubharmonic equation $\rho$ for $S$:
$S=\{\rho = 0\}$. Consider the curve
$$\gamma = \{z_j =
\gamma_j(t),\  j=1,\ldots,n, \ t\in (-\e,\e)\}\subset
S.$$
Assume that $\gamma$ is real analytic, so that locally there
exists a complex extension $\widetilde \gamma \supset \gamma$.
Then one side of $\widetilde \gamma$ lies in $\O$ if and only if
\begin{equation}\label{curve}
 \sum_j Re \frac{\partial \rho}{\partial z_j} \frac{\partial \gamma_j}{\partial t} \neq 0.
\end{equation}
Observe that condition (\ref{curve}), which depends only on $\gamma$ (when $S$ is given), is not satisfied in Example \ref{E1}. Sufficiency of (\ref{curve}) is true when $S$ is \emph{any} real hypersurface: indeed, from a geometric point of view, the condition is equivalent to the transversality of $T(\widetilde \gamma)$ and $H(S)$ (and hence $T(S)$). Pseudoconvexity is required to establish the necessity.
\er

\section{The Global Result}
In order to make the proof more transparent we first treat the
case when $\Omega$ is an unbounded convex domain with smooth
boundary $b\O$.  In the next section we will prove the main
theorem in all its generality.

\bt\label{MT} Let $M$ be a maximally complex
(connected) $(2m+1)$-real submanifold $(m \geq 1)$ of $b\O$. Assume
that $\O$ does not contain straight lines and $b\Omega=S$ satisfies the
conditions of Lemma \ref{wk}. Then there exists an $(m+1)$-complex
subvariety $W$ of $\O$, with isolated singularities, such that
$bW=M$. \et

We observe that under the hypothesis of Theorem \ref{MT}, there
exists a complex strip in a tubular neighborhood with boundary $M$
(see Corollary \ref{L4}). Moreover, since $\O$ does not contain
straight lines, we can approximate uniformly from both sides $b\O$
by strictly convex domains, see \cite{PT}. It follows that we can
find a real hyperplane $L$ such that, for any translation $L'$ of
$L$, $L'\cap \oli \O$ is a compact set. We choose an exhaustive
sequence $L_k$ of such hyperplanes, and we set $\O_k$ as the
bounded connected component of $\O\setminus L_k$. Then,
approximating from inside, we can choose a strictly convex open
subset $\O_k'\subset \O$ such that $b\O_k' \cap \O_k\subset I$,
where $I$ is the tubular neighborhood of Corollary \ref{L4}. It is
easily seen, then, that we are in the situation of the following

\bp\label{P} Let $D\Subset
B\Subset\C^n$ ($n\geq4$) be two strictly convex domains. Let
$D_+=D\cap\left\{{\sf Re} \ z_n>0\right\}$, $B_+=B\cap\left\{{\sf
Re} \ z_n>0\right\}$. Then every $(m+1)$-complex subvariety
$(m\geq2)$ with isolated singularities, $A \subset B_+ \setminus
\oli{D}_+\doteqdot C_+$, is the restriction of a complex
subvariety $\widetilde{A}$ of $B_+$ with isolated singularities.
\ep

We treat the cases $m\geq2$ and $m=1$ separately. Indeed all the
main ideas of the proof lie in the case $m\geq2$, while the case
$m=1$ simply adds technical difficulties.

\subsection{$M$ is of dimension at least $5$: $m\geq2$}
Before proving Proposition \ref{P}, we make some considerations
and we prove two lemmata that will be useful.

Let $\varphi$ be a strictly convex function\footnote{In the general case $\varphi$ will be a strongly plurisubharmonic function.} defined in a neighborhood
of $B$ such that $B=\left\{\varphi<0\right\}$. Fixing $\varepsilon>0$
small enough, $B'=\left\{\varphi<-\varepsilon\right\}$ is a strictly convex
domain of $B$ whose boundary $H$ intersects $A$ in a smooth
maximally complex submanifold $N$. A  natural way to proceed is to slice $N$ with complex
hyperplanes, in order to apply Harvey-Lawson's theorem. Each slice
of $B'$ is strictly convex, hence strongly pseudoconvex, and so
the holomorphic chain we obtain is contained in $B'$. Thus the set
made up by collecting the chains is contained in $B'$. Analyticity
of this set is the hard part of the proof.

Because of Sard's lemma, for all $ z\in D_+$, there exist a vector
$v$ arbitrarily close to $\der z_n$, and $k\in\C$ such that $z\in
v_k\doteqdot v^{\perp}+k$ and $A_k\doteqdot v_k\cap N$ is
transversal and compact, and thus smooth.

In a neighborhood of each fixed $z_0\in D_+$, the same vector $v$
realizes the transversality condition. Hence we should now fix our
attention to a neighborhood of the form
$\widehat{U}\doteqdot\bigcup_{k\in U}v_k\cap B_+$, where $v_{k_0}$ is the vector corresponding to $z_0$ and $U\subset\C$ a neighborhood of $k_0$.

Let
$\pi:\widehat{U}\to\C^{m}$ be a generic projection: we use
$(w',w)$ as holomorphic coordinates on
$v_{k_{0}}=\C^m\times\C^{n-m-1}$ (and also for $k$ near to $k_0$).
Let $V_k= \C^m \setminus \pi(A_k)$, and $V=\cap_k V_k$.

Since $A_{k_{0}}$ has a local extension (given by $v_{k_{0}}\cap
A$), it is maximally complex and so, by Harvey-Lawson's theorem,
there is a holomorphic chain $\widetilde{A}_{k_{0}}$ with
$b\widetilde{A}_{k_{0}}=A_{k_{0}}$, which extends holomorphically
$A_{k_{0}}$.

Our goal is to show that $\widetilde A_U=\cup_k\widetilde A_k$ is
analytic in $\pi^{-1}(V)$. From this, it will follow that
$\widetilde A_U$ is an analytic subvariety of $\widehat{U}$, $\pi$
being a generic projection.

Following an idea of Zaitsev, for $k\in U$,
$w'\in\C^{m}\setminus\pi(A_k)$ and $\alpha\in\N^{n-m-1}$, we
define $$ I^\alpha(w',k)\ \doteqdot\ \int_{(\eta',\eta)\in
A_k}\eta^\alpha\omega_{BM}(\eta'-w'), $$ $\omega_{BM}$ being the
Bochner-Martinelli kernel.

\bl[Zaitsev] \label{Zaitsev} Let $F(w',k)$ be the multiple-valued
function which represents $\widetilde A_k$ on
$\C^{m}\setminus\pi(A_k)$; then, if we denote by
$P^\alpha(F(w',k))$ the sum of the $\alpha^\emph{th}$ powers of the
values of $F(w',k)$, the following holds: $$ P^\alpha(F(w',k)) =
I^\alpha(w',k). $$ In particular, $F(w',k)$ is finite. \el \begin{proof}
Let $V_0$ be the unbounded component of $V_k$ (where, of course,
$P^\alpha(F(w',k)) = 0$). It is easy to show, following \cite
{HL}, that on $V_0$ also $I^\alpha(F(w',k)) = 0$: in fact, if $w'$
is far enough from $\pi(A_k)$, then $\beta = \eta^\alpha
\omega_{BM}(\eta' - w')$ is a regular $(m,m-1)$-form on some Stein
neighborhood $O$ of $A_k$. So, since in $O$ there exists $\gamma$
such that $\oli\partial\gamma = \beta$, we may write in the
language of currents $$[A_k](\beta) =
[A_k]_{m,m-1}(\oli\partial\gamma) =
\oli\partial[A_k]_{m,m-1}(\gamma) = 0.$$

In fact, since $A_k$ is maximally complex, $[A_k]=[A_k]_{m,m-1} +
[A_k]_{m-1,m}$ and $\oli\partial [A_k]_{m,m-1} = 0$, see
\cite{HL}. Moreover, since $[A_k](\beta)$ is analytic in the
variable $w'$, $[A_k](\beta)=0$ for all $w'\in V_0$.

To conclude our proof, we just need to show that the \lq\lq
jumps\rq\rq\ of the functions $P^\alpha(F(w',k))$ and
$I^\alpha(w',k)$ across the regular part of the common boundary of
two components of $V_k$ are the same.

So, let $z'\in\pi(A_k)$ be a regular point in the common boundary
of $V_1$ and $V_2$. Locally in a neighborhood of $z'$, we can
write $\widetilde A_k$ as a finite union of graphs of holomorphic
functions, whose boundaries $A_k^i$ are either in $A_k$ or empty.
In the first case, the $A_k^i$ are $CR$ graphs over $\pi(A_k)$ in
the neighborhood of $z'$. We may thus consider the jump $j_i$ of
$P^\alpha(F(w',k))$ due to a single function. We remark that the
jump for a function $f$ is $j_i=f(z')^\alpha$. The total jump will
be the sum of them.

To deal with the jump of $I^\alpha(w',k)$ across $z'$, we split
the integration set in the sets $A_k^i$ (thus obtaining the
integrals $I_i^\alpha$) and $A_k\setminus\cup_i A_k^i$
($I_0^\alpha$). Thanks to Plemelj's formulas (see~\cite{HL}) the
jumps of $I_i^\alpha$ are precisely $j_i$. Moreover, since the
form $\eta^\alpha \omega_{BM}(\eta' - z')$ is $C^\infty$
in a neighborhood of $A_k\setminus\cup_i A_k^i$, the jump of
$I_0^\alpha$ is $0$. So $P^\alpha(F(w',k))=I^\alpha(w',k)$.
\end{proof}

\br Lemma \ref{Zaitsev} implies, in particular, that the functions
$P^\alpha(F(w',k))$ are continuous in $k$. Indeed, they are
represented as integrals of a fixed form over submanifolds $A_k$
which vary continuously with the parameter $k$.\er

The functions $P^\alpha(F(w',k))$ and the holomorphic chain
$\widetilde{A}_{k_{0}}$ uniquely determine each other and so,
proving that the union over $k$ of the $\widetilde{A}_{k}$ is an
analytic set is equivalent to proving that the functions
$P^\alpha(F(w',k))$ are holomorphic in the variable $k\in
U\subset\C$.

\bl $P^\alpha(F(w',k))$ is holomorphic in the variable $k\in
U\subset\C$, for each $\alpha\in\N^{n-m-1}$. \el \begin{proof} The proof
is very similar to the one of Lewy's main lemma in~\cite{Le}. Let
us fix a point $\left(w',\underline{k}\right)$ such that $w'\notin
A_{\underline{k}}$ (this condition remains true for $k\in
B_\epsilon(\underline{k})$). Consider as domain of $P^\alpha(F)$
the set $\left\{w'\right\}\times B_\epsilon(\underline{k})$. In
view of Morera's theorem, we need to prove that for any simple
curve $\gamma\subset B_\epsilon(\underline{k})$, $$ \int_\gamma
P^\alpha(F(w',k))dk\ =\ 0. $$ Let $\Gamma\subset
B_\epsilon(\underline{k})$ be an open set such that
$b\Gamma=\gamma$. By $\gamma\ast A_k$ ($\Gamma\ast A_k$) we mean
the union of $A_k$ along $\gamma$ (along $\Gamma$). Note that
these sets are submanifolds of $N$ ($\Gamma\ast A_k$ is an open
subset) and $b(\Gamma\ast A_k)=\gamma\ast A_k$. By
Lemma~\ref{Zaitsev} and Stoke's theorem
\begin{eqnarray}
    \nonumber\int_\gamma P^\alpha(F(w',k))dk\ &=& \int_\gamma I^\alpha(w',k)dk\ =\\
    \nonumber&=&\ \int_\gamma\left(\int_{(\eta',\eta)\in A_k} \eta^\alpha\omega_{BM}(\eta'-w')\right)dk\ =\\
    \nonumber&=&\ \iint_{\gamma\ast A_k}\eta^\alpha\omega_{BM}(\eta'-w')\wedge dk\ =\\
    \nonumber&=&\ \iint_{\Gamma\ast A_k}d\left(\eta^\alpha\omega_{BM}(\eta'-w')\wedge dk\right)\ =\\
    \nonumber&=&\ \iint_{\Gamma\ast A_k}d \eta^\alpha\wedge\omega_{BM}(\eta'-w')\wedge dk\ =\\
    \nonumber&=& 0.
\end{eqnarray}
The last equality follows from the fact that since $\eta^\alpha$ is holomorphic, only holomorphic differentials appear in $d\eta^\alpha$. Since all the holomorphic differentials supported by $\Gamma\ast A_k$ already appear in $\omega_{BM}(\eta'-w')\wedge dk$, the integral is zero.
\end{proof}

We may now prove Proposition~\ref{P}.\vspace{0,3cm}

\begin{proof} \textbf{(Proposition~\ref{P}, $m\geq2$)} Up to this point we
have extended the complex manifold $A$ to an analytic set
$$\widetilde{A}_U\doteqdot A\cup\bigcup_{k\in
U}\widetilde{A}_k\subset V_U\doteqdot C_+\cup\bigcup_{k\in
U}\left(v_k\cap B_+\right).$$ The open sets $V_U$ are an open
covering of $B_+$.

Moreover the open sets
$\omega_U\doteqdot\bigcup_{k\in U}(v_k\cap B_+)$ are an open
covering of each compact set $K_\delta\doteqdot \overline
B'\cap\left\{{\sf Re}\, z_n\geq\delta\right\}$. Hence there
exist $\omega_1,\dots,\omega_l$ which cover $K_\delta$ and such
that $\omega_i\cap\omega_{i+1}\cap C_+\neq\emptyset$, for $
i=1,\dots,l-1$ and therefore there exists a countable open cover
$\left\{\omega_i\right\}_{i\in\N}$ of $\overline
B'\cap B_+$ such that, for all $i\in\N$,
$\omega_i\cap\omega_{i+1}\cap C_+\neq\emptyset$.

So we may
extend $A$ to $C_+\cup\omega_1$ by proceeding as above.

 Suppose
now that we have extended $A$ to $C^i\doteqdot
C_+\cup\bigcup^{i}_{j=1}\omega_j$ with an analytic set $A_i$. On
the non-empty intersection $C^i\cap\omega_{i+1}\cap C_+$ $A_i$ and
the extension $\widetilde{A}_{i+1}$ of $A$ to
$C_+\cup\omega_{i+1}$ coincide (as they both coincide with $A$),
hence by analicity they coincide everywhere. Consequently we may
extend $A$ to $C^{i+1}$ by $A_{i+1}\doteqdot
A_i\cup\widetilde{A}_{i+1}$. It follows that, defining $$\widetilde{A}\
\doteqdot \ A\cup\bigcup_{j\in \N}A_j,$$
$\widetilde{A}$ is the
desired extension of $A$ to $B_+$. In order to conclude the proof
we have to show that $\widetilde A$ has isolated singularities.
Let ${\sf Sing} \ (\widetilde A)\subset B'_+$ be the singular
locus of $\widetilde A$.

Recall that $\varphi$ is a strictly convex defining function for $B$. Let us consider the family
$$(\phi_\lambda\ =\ \lambda\varphi+(1-\lambda){\sf Re}\, z_n)_{\lambda\in[0,1]}$$
of strictly convex functions. For $\lambda$ near to $1$,
$\left\{\phi_\lambda=0\right\}$ does not intersect the singular
locus ${\sf Sing} \ (\widetilde A)$. Let $\oli \lambda$ be the
biggest value of $\lambda$ for which $\{\phi_\lambda=0\}\cap {\sf
Sing} \ (\widetilde A)\neq\emptyset$. Then $$\left\{\phi_{\oli
\lambda}<0\right\}\cap B_+\subset B_+$$ is a Stein domain in whose
closure the analytic set ${\sf Sing} \ (\widetilde A)$ is
contained, touching the boundary in a point of strict
convexity. So, by Kontinuit\"atsatz,
$$\{\phi_{\oli
\lambda}=0\}\cap {\sf Sing} \ (\widetilde A)$$
is a set of isolated
points in ${\sf Sing} \ (\widetilde A)$. By repeating the
argument, we conclude that ${\sf Sing} \ (\widetilde A)$ is made
up by isolated points.
\end{proof}

\begin{proof} \textbf{(Theorem~\ref{MT}, $m\geq2$)} Thanks to
Corollary~\ref{L4}, we have a regular submanifold $W_1$ of a
tubular neighborhood $I$, with boundary $M$.

Suppose $0\in M$. The real hyperplanes $H_k\doteqdot T_0(S)+k$, $k\in\R$, intesect $S$ in a compact set. If the intersection is non-empty, $\oli\O$ is divided in two sets. Let $\O_k$ be the compact one. We can choose a sequence $H_{k_n}$ such that $\O_{k_n}$ is an exaustive sequence for $\oli\O$.

We apply proposition~\ref{P} with $B_+=\O_{k_n}$, $C_+=I\cap\O_{k_n}$, and $A=W_1\cap\O_{k_n}$, to obtain an extension of $W_1$ in $\O_{k_n}$. Since, by the identity principle, two such extensions coincide in $\O_{k_{\min\left\{n,m\right\}}}$, their union is the desired submanifold $W$.
\end{proof}

\subsection{$M$ is of dimension $3$: $m=1$}
We prove now the statement of Proposition~\ref{P} for $m=1$.

Our first step is to show that when we slice transversally $N$ with complex
hyperplanes, we obtain $1$-real submanifolds which satisfy the moments condition.

Again, we fix our attention to a neighborhood of the form
$\widehat{U}\doteqdot\bigcup_{k\in U}v_k\cap B_+$. In $\widehat{U}$,
with coordinates $w_1,\ldots,w_{n-1},k$, we choose an arbitrary
holomorphic $(1,0)$-form which is constant with respect to $k$.

\bl\label{omeol}
The function
$$\Phi_\omega(k)\ =\ \int_{A_k}\omega$$
is holomorphic in $U$.
\el
\begin{proof}
We use again Morera's theorem. We need to prove that for any simple curve $\gamma\subset U$, $\gamma=b\Gamma$,
$$
\int_\gamma \Phi_\omega(k)dk\ =\ 0.
$$
Applying Stoke's theorem, we have
\begin{eqnarray}
    \nonumber\int_\gamma \Phi_\omega(k)dk\ &=& \int_\gamma\left(\int_{A_k}\omega\right)dk\ =\\
    \nonumber &=& \iint_{\gamma\ast A_k}\omega\wedge dk\ =\\
    \nonumber &=& \iint_{\Gamma\ast A_k}d(\omega\wedge dk)\ =\\
    \nonumber &=& \iint_{\Gamma\ast A_k}\partial\omega\wedge dk\ =\\
    \nonumber &=& \ 0.
\end{eqnarray}
The last equality is due to the fact that $\Gamma\ast A_k\subset N$ is maximally complex and thus supports only $(2,1)$ and $(1,2)$-forms, while $\partial\omega\wedge dk$ is a $(3,0)$-form.
\end{proof}

Now we can prove Proposition~\ref{P} and Theorem~\ref{MT} also when $m=1$.

We can find a countable covering of $B_+$ made of open subsets $\omega_i=\widehat{U}_i\cap B_+$ in such a way that:
\begin{enumerate}
    \item $\omega_0\subset C_+$;
    \item if
    $$B_l\ =\ \bigcup_{i=1}^l\omega_i,$$
    then $\omega_{l+1}\cap B_l\supset v_{l+1}\cap B_+$, where $v_{l+1}$ is a complex hyperplane in $\widehat{U}_{l+1}$.
\end{enumerate}
Now, suppose we have already found $\widetilde A_l$ that extends
$A$ on $B_l$ (observe that in $B_0=\omega_0$, $\widetilde A_0 =A$).
To conclude the proof we have to find $\widetilde A_{l+1}$
extending $A$ on $B_{l+1}$.

Each slice of $N$ in $B_l$ is maximally complex, and so are $v_{l+1}\cap N$ and $v_\epsilon\cap N$, for $v_\epsilon\subset\omega_{l+1}$ sufficiently near to $v_{l+1}$ (because they are in $B_l$ as well).

Now we use Lemma~\ref{omeol} with $\widehat U=\widehat U_{l+1}$. What
we have just observed implies that, for all holomorphic
$(1,0)$-form $\eta$, $\Phi_\eta(k)$ vanishes in an open subset of
$U$ and so is identically zero on $U$. This implies that all
slices in $\omega_{l+1}$ are maximally complex. Again we may apply
Harvey-Lawson's theorem slice by slice and conclude by the methods
of Proposition~\ref{P}.

\subsection{$M$ is of dimension $1$: $m=0$}
We have already observed that if $M$ is one-dimensional the local extension inside $\O$ may not exist (see Example~\ref{E1}). Even though there is a local strip in which we have an extension, the methods used to prove Proposition~\ref{P} do not work, since the transversal slices $M$ are either empty or isolated points. Indeed, as the following example shows, that extension result does not hold for $m=0$.

\bex \label{E2}Using the notation of Proposition~\ref{P}, in
$\C^2$ let $B$ and $D$ be the balls
$$B=\left\{|z_1|^2+|z_2|^2<c\right\},\ \ \
D=\left\{|z_1|^2+|z_2|^2<\e\right\},\ \ \ c>\e>2.$$

Consider the connected irreducible analytic set of codimension one
$$A=\{(z_1,z_2)\in B_+\ :\ z_1z_2=1\}$$ and its restriction $A_C$
to $C_+$. If $A_C$ has two connected components, $A_1$ and $A_2$,
when we try to extend $A_1$ (analytic set of codimension one on
$C_+$) to $B_+$, its restriction to $C_+$ will contain also $A_2$.
So $A_1$ is an analytic set of codimension one on $C_+$ that does
not extend on $B_+$.

So, let us prove that $A_C$ has indeed two connected components. A point of $A$ (of $A_C$) can be written as $z_1=\rho e^{i\theta}$, $z_2=\frac{1}{\rho} e^{-i\theta}$, with $\rho\in\R^+$ and $\theta\in\left(-\frac\pi2,\frac\pi2\right)$. Hence, points in $A_C$ satisfy
$$
2<\varepsilon<\rho^2+\frac1{\rho^2}<c\ \Rightarrow\ 2<\sqrt{\varepsilon+2}<\rho+\frac1\rho<\sqrt{c+2} .
$$
Since $f(\rho)=\rho+1/\rho$ is monotone decreasing up to $\rho=1$ (where $f(1)=2$), and then monotone increasing, there exist $a$ and $b$ such that the inequalities are satisfied when $a<\rho<b<1$, or when $1<1/b<\rho<1/a$. $A_C$ is thus the union of the two disjoint open sets
$$
\xymatrix{A_1=\left\{ \left(\rho e^{i\theta},\frac1\rho e^{-i\theta}\right)\in \C^2\ \Big|\ a<\rho<b,\  -\frac\pi2<\theta<\frac\pi2\right\};\\
A_2=\left\{ \left(\rho e^{i\theta},\frac1\rho e^{-i\theta}\right)\in \C^2\ \Big|\ a<\frac1\rho<b,\ -\frac\pi2<\theta<\frac\pi2\right\}.}$$

\eex
\section{Extension to Pseudoconvex Domains}
We may now prove

\bt Let $\O$ be an unbounded domain in $\C^n$ $(n\geq
3)$ with smooth boundary $b\O$ and $M$ be a maximally complex
closed $(2m+1)$-real submanifold $(m \geq 1)$ of $b\O$. Assume
that
\begin{enumerate}
\item [\emph{(i)}] $b\O$ is weakly pseudoconvex and the Levi form $\mathcal
L(b\O)$ has at least $n-m$ positive eigenvalues at every point of
$M$;
\item[\emph{(ii)}] $M$ satisfies condition $(\star)$.
\end{enumerate} Then there exists a unique
$(m+1)$-complex analytic subvariety $W$ of $\O$, such that $bW =
M$. Moreover the singular
locus of $W$ is discrete and the closure of $W$ in $\oli \O
\setminus {\sf Sing} \ W$ is a smooth submanifold with boundary
$M$. \et

\begin{proof}
Assume, for the moment, that condition ($\star$) is
replaced by the stronger condition
\begin{itemize}\item[] ${\oli \O}^\infty \cap \S_0 =
\emptyset$ where ${\oli \O^\infty}$ denotes the projective closure
of $\O$.\end{itemize}

The only thing we have to show in order to conclude the proof (by
using the methods of the previous section) is that, up to a
holomorphic change of coordinates and a holomorphic embedding $V:\C^n\rightarrow \C^N$, we can choose a sequence of real hyperplanes
$H_k\subset\C^N$, $k\in\N$, which are exhaustive in the following
sense:
\begin{itemize}
\item[1.] $H_k\cap V(S)$ is compact, for all $k\in\N$;
\item[2.] one of the two halfspaces in which $H_k$ divides $\C^N$, say $H_k^+$, intersects $V(\O)$ in a relatively compact set;
\item[3.] $\cup_k (H_k^+\cap V(\O))=V(\O)$.
\end{itemize}
The arguments of Proposition \ref{P}, indeed ---excluded the proof
that the singularities are isolated--- depend only on the fact that
we can cut $M$ by complex hyperplanes, obtaining compact maximally
complex submanifolds. Once we have found $W'\subset V(\C^n)$ ($W'$
is in fact contained in $V(\C^n)$ by analytic continuation, since
it has to coincide with the strip in a neighborhood of $M$), we
set $W= V^{-1}(W')$. Observe that the hypersurfaces $V^{-1}(H_k)$
are an exhaustive sequence for $\O$; let $\O_k$ be correspondent
sequence of relatively compact subsets. Since $\O$ is a domain of
holomorphy, for each $k$ we can choose a strongly pseudoconvex
open subset $\O_k'\subset \O$ such that $b\O_k' \cap \O_k \subset
I$, where $I$ is the tubular neighborhood found in Corollary \ref{L4}. So, in each $\O_k$ we can suppose that we deal with a
strongly pseudoconvex open set, and thus the proof of the fact
that the singularities are isolated is the same as in Proposition
\ref{P}.

Following~\cite{L2} we divide the proof in two steps.

\emph{Step 1}. $P$ linear. We consider
$\oli\O\subset\C\P^n=\C^n\cup\C\P^{n-1}_\infty$, which is disjoint
from $\S_0=\left\{P=0\right\}$. So we can consider new coordinates
of $\C\P^n$ in such a way that $\S_0$ is the $\C\P^{n-1}$ at
infinity. Now $\O$ is a relatively compact open set of
$(\C^n)'=\C\P^n\setminus \S_0$, and
$H_\infty=\C\P^{n-1}_\infty\cap(\C^n)'$ is a complex hyperplane
containing the boundary of $S$. Let $H^\R_\infty\supset H_\infty$
be a real hyperplane. The intersection between $S$ and a
translated of $H^\R_\infty$ is either empty or compact. For all
$z\in\O$, there exist a real hyperplane $H^\R_\infty\not\ni z$,
intersecting $\O$, and a small translated $H_{\e_z}$ such that
$z\in H_{\e_z}^+$. Since $\O=\cup_z (H_{\e_z}^+\cap\O)$, and $\O$
is a countable union of compact sets, we may choose an exhaustive
sequence $H_k$.

\emph{Step 2}. $P$ generic. We use the Veronese map $v$ to embed
$\C\P^n$ in a suitable $\C\P^N$ in such a way that $v(\S_0)=L_0\cap
v(\C\P^n)$, where $L_0$ is a linear subspace. The Veronese map $v$
is defined as follows: let $d$ be the degree of $P$, and let $$N\
=\ {{n+d}\choose{d}}-1. $$ Then $v$ is defined by $$v(z)\ =\
v[z_0:\ldots:z_n]\ =\ [\ldots:w_I:\ldots]_{|I|=d},$$ where
$w_I=z^I$. If $P=\sum_{|I|=d}\alpha_I z^I$, then $v(\S_0)=L_0\cap
v(\C\P^n)$, where $$L_0=\left\{\sum_{|I|=d} \alpha_I w_I = 0
\right\}.$$ Again we can change the coordinates so that $L_0$ is
the $\C\P^{N-1}$ at infinity. We may now find the exhaustive
sequence $H_k$ as in Step 1.

This achieves the proof in the case when ${\oli\O}^\infty \cap \S_0
= \emptyset$. The general case is now an easy consequence.

Indeed, since $\C\P^n\setminus \S_0$ is Stein, there is a strictly
plurisubharmonic exhaustion function $\psi$. The sets $$\O_c\ =\
\left\{\psi<c\right\}$$ are an exhaustive strongly pseudoconvex
family for $\C\P^n\setminus \S_0$. Thus in view of ($\star$) there exists $\oli c$ such
that $\oli M\subset \O_{\oli c}$. $\O'\doteqdot\O\cap\O_{\oli c}$, up to a
regularization of the boundary, is a strongly pseudoconvex open
set verifying ($\star$) in whose boundary lies $M$, and thus $M$
can be extended thanks to what has already been proved.
\end{proof}

\section*{Acknowledgments}
This research was partially supported by the MIUR project \lq\lq Geometric properties of real and complex manifolds\rq\rq.

We wish to thank Giuseppe Tomassini, whose kind help made this work possible.

Useful remarks by the referee helped us to make clearer proofs in the first part of the article and to correct various misprints.

\end{document}